\documentclass{article}

\usepackage{amsthm}
\usepackage{amsmath}
\usepackage{amssymb}
\usepackage{mathrsfs}
\usepackage{dsfont}
\usepackage[all]{xy}
\usepackage{tikz-cd}
\usepackage{enumitem}
\usepackage{stmaryrd}
\usepackage{hyperref}
\usepackage{newclude}

\theoremstyle{definition}
\newtheorem{Definition}{Definition}[subsection]

\theoremstyle{plain}
\newtheorem{Theorem}[Definition]{Theorem}

\theoremstyle{plain}

\theoremstyle{plain}
\newtheorem{Proposition}[Definition]{Proposition}

\theoremstyle{plain}

\theoremstyle{plain}
\newtheorem{Lemma}[Definition]{Lemma}

\theoremstyle{plain}

\theoremstyle{plain}

\theoremstyle{plain}

\theoremstyle{plain}

\theoremstyle{plain}

\theoremstyle{definition}
\newtheorem{Example}[Definition]{Example}

\theoremstyle{definition}
\newtheorem{Notation}[Definition]{Notation}

\theoremstyle{remark}
\newtheorem{Remark}[Definition]{Remark}

\theoremstyle{plain}
\newcommand{\thistheoremname}{}
\newtheorem*{genericthm*}{\thistheoremname}
\newenvironment{namedthm*}[1]
  {\renewcommand{\thistheoremname}{#1}%
   \begin{genericthm*}}
  {\end{genericthm*}}

\author{Thibault D. D\'ecoppet}
\title{Multifusion Categories and Finite Semisimple 2-Categories}

\begin{document}

\bibliographystyle{alpha}




    \maketitle
    \hspace{1cm}
    \begin{abstract}
        We give a 3-universal property for the Karoubi envelope of a 2-category. Using this, we show that the 3-categories of finite semisimple 2-categories (as introduced in \cite{DR}) and of multifusion categories are equivalent.
    \end{abstract}

\section*{Introduction}

After Lurie's sketch of proof of the Cobordism Hypothesis of Baez and Dolan (see \cite{L}), there began a quest for interesting fully-dualizable objects in symmetric monoidal higher categories. The most well-known examples are perhaps given by fusion categories (over an algebraically closed field of characteristic zero), which are fully-dualizable objects of the symmetric monoidal 3-category $\mathbf{Mult}$ of (multi)fusion categories, bimodules categories, bimodule functors, and bimodule natural transformations with monoidal structure given by the Deligne tensor product (see \cite{DSPS13}). As theorem 1.4.8 of \cite{DR} shows, finite semisimple 2-categories essentially correspond to multifusion categories. Further, as finite semisimple 2-categories can be assembled into a 3-category $\mathbf{FSS2C}$, it is natural to ask how it compares with the 3-category of multifusion categories. We show that the canonical 3-functor sending a multifusion categories to its 2-category of finite semisimple module categories induces an equivalence $\mathbf{Mult}\simeq\mathbf{FSS2C}$. In subsequent work, we will show how to define a product on the 3-category of finite semisimple 2-categories, and explain how it relates to the Deligne tensor product of multifusion categories under the above equivalence.

We now describe the content of the different sections of the present article. Section \ref{sec:2condensation} begins by recalling the definitions of a 2-condensation monad and that of a 2-condensation given in \cite{GJF}. As explained there, 2-condensation monads categorify the notion of an idempotent in a 1-category. In particular, we recall what it means for a 2-condensation monad to split. In \cite{GJF}, it is shown that locally idempotent complete 2-categories admit a Karoubi envelope, i.e. a 2-functor to an associated 2-category in which every 2-condensation monad splits. We state the 3-universal property of this construction, which will be proven in appendix \ref{sec:completion2categories}.

In section \ref{sec:fusion}, we present a weak form of the definitions of (finite) semisimple 2-categories given in \cite{DR}: We don't impose any strictness hypothesis on the underlying 2-category and use 2-condensation monads instead of separable monads. We show that every finite semisimple 2-category (as defined in the present text) is equivalent to a finite semisimple 2-category as defined by Douglas-Reutter. Then, we proceed to prove our main theorem.

\begin{namedthm*}{Theorem \ref{thm:comparaison}}
There is a linear 3-functor $$\mathbf{Mod}:\mathbf{Mult}\rightarrow \mathbf{FSS2C}$$ that sends a multifusion category to its 2-category of finite semisimple right module categories. Further, this 3-functor is an equivalence.
\end{namedthm*}

\noindent Finally, we define connected finite semisimple 2-categories, and explain what the corresponding property at the level of multifusion categories is.

\section{2-Condensation Monads and 2-Condensations}\label{sec:2condensation}

\subsection{Definitions}

The notion of an $n$-condensation monad in an $n$-category was introduced in \cite{GJF}. These objects categorify the notion of an idempotent in a (1-)category to an $n$-category. In particular, a 1-condensation monad is an idempotent in a (1-)category. Below, we recall the definition in the case $n=2$.

\begin{Definition}
Let $\mathfrak{C}$ be a 2-category. A 2-condensation monad in $\mathfrak{C}$ consists of:
\begin{enumerate}
    \item An object $A$;
    \item A 1-morphism $e:A\rightarrow A$;
    \item Two 2-morphisms $\mu: e\circ e\Rightarrow e$ and $\delta:e\Rightarrow e\circ e$,
\end{enumerate}

such that

\begin{enumerate}
    \item The pair $(e,\mu)$ is a (non-unital) associative algebra;
    \item The pair $(e,\delta)$ is a (non-counital) coassociative coalgebra;
    \item The Frobenius relations hold, i.e. $\delta$ is an $(e,e)$-bimodule map;
    \item We have: $\mu \cdot \delta = Id_e$.
\end{enumerate}
\end{Definition}

\begin{Remark}
More concisely, a 2-condensation monad on $A$ in $\mathfrak{C}$ is a non-unital non-counital special Frobenius algebra in the monoidal category $End_{\mathfrak{C}}(A)$.
\end{Remark}

Generalizing the notion of a split surjection in a (1-)category, \cite{GJF} also introduced the notion of an $n$-condensation in an $n$-category. The case $n=2$ is reviewed in the definition below.

\begin{Definition}
Let $\mathfrak{C}$ be a 2-category. A $2$-condensation in $\mathfrak{C}$ consists of:
\begin{enumerate}
    \item Two objects $A$ and $B$;
    \item Two 1-morphism $f:A\rightarrow B$ and $g:B\rightarrow A$;
    \item Two 2-morphisms $\phi: f\circ g\Rightarrow Id_B$ and $\gamma:Id_B\Rightarrow f\circ g$
\end{enumerate}

such that

\begin{enumerate}
    \item We have: $\phi \cdot \gamma = Id_{Id_B}$.
\end{enumerate}
\end{Definition}

Given a 2-condensation as above, observe that the induced structure on the object $A$ alone is precisely that of a 2-condensation monad.

\begin{Definition}
Given a 2-condensation monad $(A,e,\mu,\delta)$ in $\mathfrak{C}$, a splitting, or, more precisely, an extension of this 2-condensation monad to a 2-condensation, is a 2-condensation $(A,B,f,g,\phi,\gamma)$ together with a 2-isomorphism $\theta:g\circ f\cong e$ such that
$$\mu = \theta\cdot (g\circ \phi\circ f)\cdot(\theta^{-1} \circ \theta^{-1}),$$
$$\delta = (\theta \circ \theta)\cdot (g\circ \gamma\circ f)\cdot\theta^{-1}.$$
\end{Definition}

\begin{Remark}
Extensions of a 2-condensation monad to a 2-condensation are preserved by 2-functors (upon insertion of the relevant coherence 2-isomorphisms).
\end{Remark}

Categorifying the notion of being idempotent complete \cite{GJF} makes the following definition.

\begin{Definition}
A 2-category $\mathfrak{C}$ has all condensates if it is locally idempotent complete and every 2-condensation monad can be extended to a 2-condensation.
\end{Definition}

\begin{Remark}\label{rem:uniquenesssplitting2condensationmonads}
Theorem 2.3.2 of \cite{GJF} states that if $\mathfrak{C}$ is a locally idempotent complete 2-category, then the 2-category of extensions of any 2-condensation monad to a 2-condensation is either empty or contractible. (This is the categorified version of the statement that if an idempotent in a (1-)category splits, it does so uniquely up to unique isomorphism.)
\end{Remark}

\subsection{The Karoubi Envelope of a 2-Category}

Given a locally idempotent complete 2-category $\mathfrak{C}$, \cite{GJF} constructs a Karoubi envelope for $\mathfrak{C}$, i.e. a 2-category $Kar(\mathfrak{C})$ that has all condensates and comes with a canonical fully faithful inclusion $\iota_{\mathfrak{C}}:\mathfrak{C}\hookrightarrow Kar(\mathfrak{C})$. This is given by the 2-category whose objects are the 2-condensation monads in $\mathfrak{C}$, 1-morphism are bimodules (see definition \ref{def:2condensationbimodule}), and 2-morphisms are bimodule maps. For our purposes, it will be important to have a universal characterization for $Kar(\mathfrak{C})$ and $\iota_{\mathfrak{C}}$.

\begin{Definition} \label{def:KarEnv}
Let $\mathfrak{C}$ be a locally idempotent complete 2-category. A Karoubi envelope is a 2-functor $\iota_{\mathfrak{C}}:\mathfrak{C}\rightarrow Kar(\mathfrak{C})$ satisfying the following 3-universal property:
\begin{enumerate}
\setcounter{enumi}{-1}
\item The locally idempotent complete 2-category $Kar(\mathfrak{C})$ has all condensates.

\item For any locally idempotent complete 2-category $\mathfrak{A}$ that has all condensates and every 2-functor $F:\mathfrak{C}\rightarrow \mathfrak{A}$, there exists a 2-functor $F':Kar(\mathfrak{C})\rightarrow\mathfrak{A}$ and a 2-natural equivalence $u:F'\circ\iota_{\mathfrak{C}}\Rightarrow F$.

\item For every 2-functors $G,H:Kar(\mathfrak{C})\rightarrow \mathfrak{A}$ and 2-natural transformation $t:G\circ\iota_{\mathfrak{C}}\Rightarrow H\circ\iota_{\mathfrak{C}}$, there exists a 2-natural transformation $t':G\Rightarrow H$ and an invertible modification $\lambda:t'\circ \iota_{\mathfrak{C}}\Rrightarrow t$.

\item Furthermore, for any 2-natural transformations $r,s:G\Rightarrow H$ and modification $\lambda:r \circ\iota_{\mathfrak{C}}\Rrightarrow s \circ\iota_{\mathfrak{C}}$, there exists a unique modification $\lambda':r \Rrightarrow s$ such that $\lambda'\circ \iota_{\mathfrak{C}}=\lambda$.
\end{enumerate}
\end{Definition}

\begin{Remark}
By its very definition, a Karoubi envelope, if it exists, is unique in the sense that the 3-category of Karoubi envelopes is a contractible 3-groupoid.
\end{Remark}

In appendix \ref{sec:completion2categories}, we shall prove the following proposition.

\begin{Proposition}\label{prop:exsitenceKaroubienvelope}
Every locally idempotent complete 2-category has a Karoubi envelope.
\end{Proposition}

In the context of $R$-linear 2-categories, where $R$ is a fixed commutative ring, the natural notion of completion is that of Cauchy completion.

\begin{Definition}
An $R$-linear 2-category is Cauchy complete if it has all condensates, is locally additive, and has all direct sums (of objects).
\end{Definition}

In the linear context, one can formulate an analogue of definition \ref{def:KarEnv}, where Karoubi envelope is replaced by Cauchy complete, and the 2-functors appearing are assumed to be $R$-linear. We prove the following result in appendix \ref{sec:completion2categories}.

\begin{Proposition}\label{prop:existenceCauchycompletion}
Let $\mathfrak{C}$ be a locally additive and locally idempotent complete $R$-linear 2-category, then there exists an $R$-linear 2-functor $\mathfrak{C}\rightarrow Cau(\mathfrak{C})$ that is a Cauchy completion.
\end{Proposition}

\section{Multifusion Categories and Finite Semisimple 2-Categories}\label{sec:fusion}

Throughout, we work over a fixed algebraically closed field $\mathds{k}$ of characteristic zero.

\subsection{Semisimple 2-Categories}

The notion of a semisimple 2-category was originally defined in section 1.4 of \cite{DR}. We now give a categorically weaker definition inspired by the existing one, and prove that they are essentially equivalent. In the mean time, in order to avoid confusion, we shall use the term strict semisimple 2-categories when we refer to the objects defined in \cite{DR}.

\begin{Definition}\label{def:semisimple2category}
A semisimple 2-category is a $\mathds{k}$-linear 2-category that is locally semisimple, has adjoints for 1-morphisms, and is Cauchy complete.
\end{Definition}

\begin{Definition}
Let $\mathfrak{C}$ be a semisimple 2-category. An object $C$ of $\mathfrak{C}$ is said to be simple if the monoidal unit of $End_{\mathfrak{C}}(C)$ is a simple object.
\end{Definition}

\begin{Definition}
A finite semisimple 2-category is a semisimple 2-category that is locally finite semisimple, and has only finitely many equivalence classes of simple objects.
\end{Definition}

\begin{Lemma}\label{lem:strictFSS2C}
Let $\mathfrak{C}$ be a strict 2-category. It is (finite) semisimple in the sense of definition \ref{def:semisimple2category} if and only if it is (finite) semisimple in the sense of \cite{DR}.
\end{Lemma}
\begin{proof}
Observe that all the defining properties of a semisimple 2-category (as in definition \ref{def:semisimple2category}) do not depend on any underlying strictness hypothesis. Now, the result follows from the fact that for a locally idempotent complete and locally additive linear 2-category, being Cauchy complete is equivalent to being additive and idempotent complete (in the sense of \cite{DR}) by theorem 3.3.3 of \cite{GJF}.
\end{proof}

\begin{Lemma}\label{lem:strictificationFSS2C}
Given a (finite) semisimple 2-category $\mathfrak{C}$ (as in definition \ref{def:semisimple2category}), there exists a linearly 2-equivalent strict (finite) semisimple 2-category.
\end{Lemma}
\begin{proof}
By a linear version of the coherence theorem for 2-categories, there exists a linear 2-equivalence $F:\mathfrak{C}\rightarrow \mathfrak{D}$, where $\mathfrak{D}$ is a strict linear 2-category. Thus, it is enough to check that all the properties in the definition of a semisimple 2-categories are invariant under equivalences of linear 2-categories. This is clear for all of them. For instance, one can see that being Cauchy complete is invariant using the 3-universal property of the Cauchy completion. To complete the proof, it is enough to appeal to lemma \ref{lem:strictFSS2C}.
\end{proof}

\begin{Remark}
In particular, all the results of \cite{DR} on strict (finite) semisimple 2-categories apply to (finite) semisimple 2-categories. Most notably, we have the following theorem.
\end{Remark}

\begin{Theorem}{\cite[Thm 1.4.8]{DR}}\label{thm:ss2cmodules}
Let $\mathcal{C}$ be a multifusion category. Then $\mathbf{Mod}(\mathcal{C})$, the 2-category of finite semisimple right $\mathcal{C}$-module categories, right $\mathcal{C}$-module functors, and right $\mathcal{C}$-module natural transformations is a finite semisimple 2-category.
\end{Theorem}

\begin{Remark}
The 2-category $\mathbf{Mod}(\mathcal{C})$ can equivalently be described as the 2-category of separable algebras in $\mathcal{C}$, bimodules, and bimodule maps. This is the content of proposition 1.3.13 of \cite{DR}. We shall alternate between these two descriptions of the finite semisimple 2-category $\mathbf{Mod}(\mathcal{C})$ depending on the situation. Furthermore, by theorem 3.3.3 of \cite{GJF}, we can also think of $\mathbf{Mod}(\mathcal{C})$ as the 2-category of 2-condensation monads in $\mathrm{B}\mathcal{C}$, the one object 2-category on $\mathcal{C}$.
\end{Remark}

\begin{Definition}\label{def:generator}
Let $\mathfrak{C}$ be a finite semisimple 2-category. We say that an object $C$ in $\mathfrak{C}$ is a generator if the canonical inclusion $$\mathrm{B}End_{\mathfrak{C}}(C)\hookrightarrow \mathfrak{C}$$ of the one object 2-category with endomorphism category $End_{\mathfrak{C}}(C)$ into $\mathfrak{C}$ is a Cauchy completion.
\end{Definition}

\begin{Example}\label{ex:generator}
Let $\mathcal{C}$ be a multifusion category. Then, $\mathcal{C}$ is a generator of $\mathbf{Mod}(\mathcal{C})$. In particular, given that $End_{\mathbf{Mod}(\mathcal{C})}(\mathcal{C})\simeq \mathcal{C}$ as monoidal categories, we have $Cau(\mathrm{B}\mathcal{C})\simeq \mathbf{Mod}(\mathcal{C})$.
\end{Example}

One of the most important results on finite semisimple 2-categories is that they are precisely the 2-categories of finite semisimple module categories over a multifusion category, as was shown in theorem 1.4.9 of \cite{DR}. We recall this result in our language now.

\begin{Lemma}\label{lem:SemsimpleKaroubi}
Every finite semisimple 2-category $\mathfrak{C}$ has a generator.
\end{Lemma}
\begin{proof}
Without loss of generality, we may assume that $\mathfrak{C}$ is strict. By theorem 1.4.9 of \cite{DR}, there exists an object $C$ of $\mathfrak{C}$ such that the inclusion $\mathrm{B}End_{\mathfrak{C}}(C)\hookrightarrow \mathfrak{C}$ is a direct sum completion followed by an idempotent completion, whence it is a Cauchy completion (by theorem 3.3.3 of \cite{GJF}).
\end{proof}

Let us examine an example of a finite semisimple 2-category in details.

\begin{Example}\label{ex:Mod(p)ss}
Let $p$ be a prime. We write $\mathbf{Vect}_{\mathbb{Z}/p\mathbb{Z}}$ for the category of $\mathbb{Z}/p\mathbb{Z}$-graded vector spaces with its standard monoidal structure. We shall denote its simple objects by $\mathds{k}_i$ with $i\in\mathbb{Z}/p\mathbb{Z}$. By theorem \ref{thm:ss2cmodules}, $\mathbf{Mod}(\mathbf{Vect}_{\mathbb{Z}/p\mathbb{Z}})$ is a finite semisimple 2-category. Using theorem 1.1 of \cite{Nat}, we see that this finite semisimple 2-category has two equivalence classes of simple objects, which we denote by $\mathcal{M}(0, triv)$ and $\mathcal{M}(\mathbb{Z}/p\mathbb{Z}, triv)$. Here $triv$ denotes the trivial cohomology class (of the appropriate cohomology group). In this case, it is easy to give explicit $\mathbf{Vect}_{\mathbb{Z}/p\mathbb{Z}}$-module categories representatives: We have equivalences $\mathcal{M}(0, triv)\simeq \mathbf{Vect}_{\mathbb{Z}/p\mathbb{Z}}$ and $\mathcal{M}(\mathbb{Z}/p\mathbb{Z}, triv)\simeq\mathbf{Vect}$. Using these representatives, one can describe the finite semisimple $Hom$-categories between the two simple objects. We obtain the schematic description below.
$$
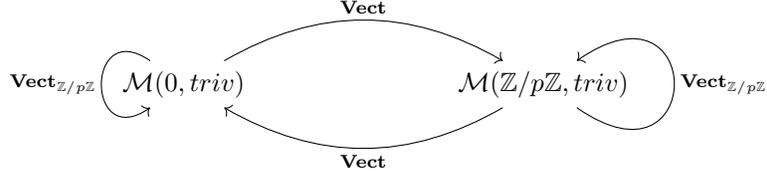
\begin{tikzcd}
\mathcal{M}(0, triv) \arrow[rrr, "\mathbf{Vect}", bend left] \arrow["\mathbf{Vect}_{\mathbb{Z}/p\mathbb{Z}}"', loop, distance=3em, in=215, out=145] &  &  & \mathcal{M}(\mathbb{Z}/p\mathbb{Z}, triv) \arrow[lll, "\mathbf{Vect}", bend left] \arrow["\mathbf{Vect}_{\mathbb{Z}/p\mathbb{Z}}"', loop, distance=6em, in=35, out=325]
\end{tikzcd}
$$
It should be noted that we have written the endomorphism category of the simple object represented by $\mathcal{M}(\mathbb{Z}/p\mathbb{Z}, triv)$ as $\mathbf{Vect}_{\mathbb{Z}/p\mathbb{Z}}$, whereas the canonical description is $Rep(\mathbb{Z}/p\mathbb{Z})$. As $\mathbb{Z}/p\mathbb{Z}$ is abelian, these two fusion categories are equivalent.

Interpreted correctly, the diagram above contains some information about the composition of 1-morphisms: Viewing $\mathbf{Vect}_{\mathbb{Z}/p\mathbb{Z}}$ as a monoidal category, and $\mathbf{Vect}$ as a $\mathbf{Vect}_{\mathbb{Z}/p\mathbb{Z}}$-bimodule category in the standard way  describes how most of the 1-morphisms compose. The remaining compositions are given by $$\begin{tabular}{ccc}$\mathbf{Vect}\times \mathbf{Vect}$&$\rightarrow$ & $\mathbf{Vect}_{\mathbb{Z}/p\mathbb{Z}}$\\$(\mathds{k},\mathds{k})$&$\mapsto$ &$\bigoplus_{i\in\mathbb{Z}/p\mathbb{Z}}\mathds{k}_i,$\end{tabular}$$
which is a $\mathbf{Vect}_{\mathbb{Z}/p\mathbb{Z}}$-balanced functor.
\end{Example}

\subsection{A Comparison of 3-Categories}

Let us begin by observing that finite semisimple 2-categories form a 3-category. Namely, it is the full sub-3-category of the linear 3-category $2Cat_{\mathds{k}}$ of $\mathds{k}$-linear 2-categories on the finite semisimple 2-categories. We denote the 3-category of finite semisimple 2-categories by $\mathbf{FSS2C}$.

In the thesis \cite{Sch}, it has been proven that fusion categories, finite semisimple bimodule categories, bimodule functors, and bimodule natural transformations form a 3-category. However, the proof supplied there does not use all the properties of fusion categories; Inspection shows that the following (marginally) stronger result was in fact proven.

\begin{Proposition}{\cite{Sch}}
Multifusion categories, finite semisimple bimodule categories, bimodule functors and bimodule natural transformations form a linear 3-category, which we denote by $\mathbf{Mult}$.
\end{Proposition}

We will now proceed to prove that these two 3-categories are equivalent, a result which was expected by the authors of \cite{DR}.

\begin{Theorem}\label{thm:comparaison}
There is a linear 3-functor $\mathbf{Mod}:\mathbf{Mult}\rightarrow \mathbf{FSS2C}$ that sends a multifusion category $\mathcal{C}$ to the associated finite semisimple 2-category of finite semisimple right $\mathcal{C}$-module categories. Moreover, this 3-functor is an equivalence.
\end{Theorem}
\begin{proof}
Let us consider the linear 3-functor $$Hom_{\mathbf{Mult}}(\mathbf{Vect},-):\mathbf{Mult}\rightarrow \mathbf{2Cat}_{\mathds{k}}.$$ By definition of the linear 3-category $\mathbf{Mult}$, this 3-functor sends a multifusion category $\mathcal{C}$ to the associated 2-category of finite semisimple $(\mathbf{Vect},\mathcal{C})$-bimodule categories. But, as a left $\mathbf{Vect}$-module category structure on a finite semisimple category is essentially unique, we may identify this 2-category with the 2-category of finite semisimple right $\mathcal{C}$-module categories. As linear 2-categories of this form are finite semisimple by theorem \ref{thm:ss2cmodules}, the 3-functor $Hom_{\mathbf{Mult}}(\mathbf{Vect},-)$ factors through the full sub-3-category $\mathbf{FSS2C}$. We denote this 3-functor by $$\mathbf{Mod}:\mathbf{Mult}\rightarrow \mathbf{FSS2C}.$$

It remains to prove that $\mathbf{Mod}$ is an equivalence of linear 3-categories. According to the definition given in \cite{GPS}, it is enough to show that $\mathbf{Mod}$ is essentially surjective and induces equivalences on $Hom$-2-categories. Lemma \ref{lem:SemsimpleKaroubi} shows that $\mathbf{Mod}$ is essentially surjective. Now, let $\mathcal{C}$ and $\mathcal{D}$ be multifusion categories. We wish to prove that the linear 2-functor $$Hom_{\mathbf{Mult}}(\mathcal{C},\mathcal{D})=\mathbf{Bimod}(\mathcal{C},\mathcal{D})\rightarrow Fun_{\mathds{k}}(\mathbf{Mod}(\mathcal{C}),\mathbf{Mod}(\mathcal{D}))$$ induced by $\mathbf{Mod}$ is a equivalence. By construction (see theorem 3.6.1 of \cite{Sch}), it sends a finite semisimple $(\mathcal{C},\mathcal{D})$-bimodule category $\mathcal{M}$ to the 2-functor $$\begin{tabular}{ccc}$\mathbf{Mod}(\mathcal{C})$&$\rightarrow$& $\mathbf{Mod}(\mathcal{D}).$\\ $\mathcal{N}$&$\mapsto$& $\mathcal{N}\boxtimes_{\mathcal{C}}\mathcal{M}$\end{tabular}$$

Firstly, we have to show that every linear 2-functor $F:\mathbf{Mod}(\mathcal{C})\rightarrow \mathbf{Mod}(\mathcal{D})$ can be obtained in this way. Observe that the image of $\mathcal{C}$ in $ \mathbf{Mod}(\mathcal{C})$ under $F$ inherits a $(\mathcal{C},\mathcal{D})$-bimodule structure from the left $\mathcal{C}$-module structure of $\mathcal{C}$. Let $\kappa:\mathrm{B}End_{\mathbf{Mod}(\mathcal{C})}(\mathcal{C})\rightarrow \mathbf{Mod}(\mathcal{C})$ be the inclusion. We obtain a strictly commutative diagram $$\xymatrix{\mathrm{B}End_{\mathbf{Mod}(\mathcal{C})}(\mathcal{C})\ar[rd]^-{F'}\ar[d]_-{\kappa}& \\
\mathbf{Mod}(\mathcal{C})\ar[r]_-F& \mathbf{Mod}(\mathcal{D})}$$ by letting $F':= F\circ\kappa.$ As $\kappa$ is a Cauchy completion (see example \ref{ex:generator}), we find that $F$ is characterized uniquely up to 2-natural equivalence by the linear 2-functor $F':\mathrm{B}End_{\mathbf{Mod}(\mathcal{C})}(\mathcal{C})\rightarrow \mathbf{Mod}(\mathcal{D})$ through the 3-universal property of the Cauchy completion (see proposition \ref{prop:existenceCauchycompletion}). Now, note that $$F'(\mathcal{C})=F(\mathcal{C})\simeq \mathcal{C}\boxtimes_{\mathcal{C}}F(\mathcal{C})$$ as $(\mathcal{C},\mathcal{D})$-bimodule categories. This means that there is a 2-natural equivalence $\kappa (-)\boxtimes_{\mathcal{C}}F(\mathcal{C})\simeq F'(-)$. Thus, by proposition \ref{prop:existenceCauchycompletion}, there is a 2-natural equivalence $F(-)\simeq (-)\boxtimes_{\mathcal{C}}F(\mathcal{C})$.

Secondly, let $\mathcal{M}_1$ and $\mathcal{M}_2$ be two $(\mathcal{C},\mathcal{D})$-bimodule categories, and $f$ be a 2-natural transformation $(-)\boxtimes_{\mathcal{C}}\mathcal{M}_1\Rightarrow (-)\boxtimes_{\mathcal{C}}\mathcal{M}_2$. Using the 3-universal property of the Cauchy completion (proposition \ref{prop:existenceCauchycompletion}), we find that $f$ is uniquely determined up to invertible modification by the functor $f_{\mathcal{C}}: \mathcal{C}\boxtimes_{\mathcal{C}}\mathcal{M}_1\rightarrow\mathcal{C}\boxtimes_{\mathcal{C}}\mathcal{M}_2$ viewed as a $(\mathcal{C},\mathcal{D})$-bimodule functor. But, we have that $$\mathcal{C}\boxtimes_{\mathcal{C}}(-):\mathbf{Bimod}(\mathcal{C},\mathcal{D})\rightarrow \mathbf{Bimod}(\mathcal{C},\mathcal{D})$$ is an equivalence of linear 2-categories, as the $(\mathcal{C},\mathcal{C})$-bimodule $\mathcal{C}$ is the identity 1-morphism on the object $\mathcal{C}$ of $\mathbf{Mult}$. Thus, we can find a $(\mathcal{C},\mathcal{D})$-bimodule functor $f':\mathcal{M}_1\rightarrow \mathcal{M}_2$ together with a bimodule natural isomorphism $\mathcal{C}\boxtimes_{\mathcal{C}}f' \cong f_{\mathcal{C}}$. Equivalently, there is an invertible modification $f_{\kappa(-)}\cong \kappa(-)\boxtimes_{\mathcal{C}}f'$. Thanks to proposition \ref{prop:existenceCauchycompletion}, this implies that there exists an invertible modification witnessing $f_{(-)} \cong (-)\boxtimes_{\mathcal{C}}f'$.

Thirdly, let $f,g:\mathcal{M}_1\rightarrow \mathcal{M}_2$ be two $(\mathcal{C},\mathcal{D})$-bimodule functors, and let $\psi$ be a modification $(-)\boxtimes_{\mathcal{C}}f\Rrightarrow (-)\boxtimes_{\mathcal{C}}g$. The 3-universal property of the Cauchy completions shows that $\psi$ determines and is uniquely determined by the $(\mathcal{C},\mathcal{D})$-bimodule natural transformation $\psi_{\mathcal{C}}:\mathcal{C}\boxtimes_{\mathcal{C}}f\Rightarrow \mathcal{C}\boxtimes_{\mathcal{C}}g$. Using again that $\mathcal{C}\boxtimes_{\mathcal{C}}(-):\mathbf{Bimod}(\mathcal{C},\mathcal{D})\rightarrow \mathbf{Bimod}(\mathcal{C},\mathcal{D})$ is an equivalence of 2-categories, we find that there exists a unique $(\mathcal{C},\mathcal{D})$-bimodule natural transformation $\psi':f\Rightarrow g$ such that $\psi_{\mathcal{C}} = \mathcal{C}\boxtimes_{\mathcal{C}}\psi'$, whence $\psi_{(-)} = (-)\boxtimes_{\mathcal{C}}\psi'$. This proves that $\mathbf{Mod}$ is fully faithful at the level of 3-morphisms, which completes the proof.
\end{proof}

\begin{Remark}
It is known that every equivalence of linear 3-categories has an inverse. We believe that the linear 3-functor $$End:\mathbf{FSS2C}\rightarrow \mathbf{Mult}$$ given on objects by sending a finite semisimple 2-category $\mathfrak{C}$ to $End_{\mathfrak{C}}(C)$, where $C$ in $\mathfrak{C}$ is an arbitrary generator, is such an inverse to $\mathbf{Mod}$.
\end{Remark}

\subsection{Connected Finite Semisimple 2-Categories}

Fusion categories are the most important examples of multifusion categories. Therefore, it seems sensible to consider the full sub-3-category of $\mathbf{Mult}$ on the fusion categories. However, this definition is problematic, as being fusion is not a Morita invariant property. Said differently, the full sub-3-category on the fusion categories is not replete.

\begin{Definition}
A multifusion category $\mathcal{C}$ is connected if it is Morita equivalent to a fusion category.
\end{Definition}

\begin{Remark}
Let $\mathcal{C}$ be a multifusion category. It is known that we can write $\mathcal{C}$ as a matrix of finite semisimple categories whose diagonal entries are non-zero fusion categories (see section 4.3 of \cite{EGNO}). It is immediate that $\mathcal{C}$ is connected if and only if all entries of this matrix are non-zero.
\end{Remark}

The folding and unfolding constructions of \cite{DR} suggest that the corresponding property for finite semisimple 2-categories is also a form of connectedness.

\begin{Definition}\label{def:connectedss2c}
Let $\mathfrak{C}$ be a finite semisimple 2-category. We say that two simple objects $A,B$ of $\mathfrak{C}$ are connected if the category $Hom_{\mathfrak{C}}(A,B)$ is non-zero. A semisimple 2-category $\mathfrak{C}$ is said to be connected if all of its simple objects are connected.
\end{Definition}

\begin{Example}
The finite semisimple 2-category of example \ref{ex:Mod(p)ss} is connected.
\end{Example}

\begin{Proposition}
The linear 3-equivalence $\mathbf{Mod}$ of theorem \ref{thm:comparaison} restricts to an equivalence between the full sub-3-categories on the connected multifusion categories and on the connected finite semisimple 2-categories.
\end{Proposition}

\begin{proof}
Observe that it is enough to show the following two things: Firstly, given a fusion category $\mathcal{C}$, then $\mathbf{Mod}(\mathcal{C})$ is connected. Secondly, if $\mathfrak{C}$ is a connected finite semisimple 2-category, then there exists a connected multifusion category $\mathcal{D}$ such that $\mathbf{Mod}(\mathcal{D})\simeq \mathfrak{C}$.

Given a fusion category $\mathcal{C}$, every finite semisimple indecomposable right $\mathcal{C}$-module category $\mathcal{M}$ admits a non-zero right $\mathcal{C}$-module functor $\mathcal{C}\rightarrow \mathcal{M}$. Namely, pick a non-zero object $M$ in $\mathcal{M}$, and consider the right module functor given by sending $C$ in $\mathcal{C}$ to $M\otimes C$. Thus, $Hom_{\mathbf{Mod}(\mathcal{C})}(\mathcal{C},\mathcal{M})$ is non-zero for every $\mathcal{M}$ as above. Now, observe that finite semisimple indecomposable $\mathcal{C}$-module categories are precisely the simple objects of $\mathbf{Mod}(\mathcal{C})$ by lemma 7.12.6 of \cite{EGNO}. As $\mathcal{C}$ is indecomposable, we are done by the existence of adjoints for $\mathcal{C}$-module functors (see \cite{DSPS14}) and proposition 1.2.19 of \cite{DR}.

Let $\mathfrak{C}$ be a connected finite semisimple 2-category. Let $A_i$ be representatives for the finitely many equivalence classes of simple objects of $\mathfrak{C}$. Let us define $$\mathcal{C}:=End_{\mathfrak{C}}(\boxplus_i A_i).$$ By definition, we have that $\mathcal{C}$ is connected. Further, the proof of theorem 1.4.9 of \cite{DR} proves that the canonical inclusion $\mathbf{Mod}(\mathcal{C})\rightarrow \mathfrak{C}$ is an equivalence.
\end{proof}

\appendix

\section{Completions}

\subsection{Completions of 1-Categories}

We recall some well-known fact about the Karoubi completions of 1-categories. Not only does this serves as a preamble to our discussion of the completions of 2-categories, but also, we will make use of some of these results later.

\begin{Definition}\label{def:1KaroubiCompletion}
Let $\mathcal{C}$ be a (1-)category. A Karoubi envelope, also called an idempotent completion, is a functor $\iota_{\mathcal{C}}:\mathcal{C}\rightarrow Kar(\mathcal{C})$ satisfying the following 2-universal property:
\begin{enumerate}
\setcounter{enumi}{-1}
\item The category $Kar(\mathcal{C})$ is idempotent complete.
\item For every idempotent complete category $\mathcal{A}$, and every functor $F:\mathcal{C}\rightarrow\mathcal{A}$, there exists a functor $F':Kar(\mathcal{C})\rightarrow\mathcal{A}$ and a natural isomorphism $t:F'\circ\iota_{\mathcal{C}}\Rightarrow F$.
\item For every functors $G,H:Kar(\mathcal{C})\rightarrow \mathcal{A}$ and natural transformation $r:G\circ\iota_{\mathcal{C}}\Rightarrow H\circ\iota_{\mathcal{C}}$, there exists a unique natural transformation $r':G\Rightarrow H$ such that $ r'\circ\iota_{\mathcal{C}}= r$.
\end{enumerate}
\end{Definition}

\begin{Remark}\label{rem:2univbiuniv}
In the terminology of \cite[Def. 9.4.]{Fio}, the 1-morphism $\iota_{\mathcal{C}}$ is biuniversal with respect to the forgetful 2-functor $U:\mathbf{2Cat}^{ic}\rightarrow \mathbf{2Cat}$ from the 2-category of (small) idempotent complete categories to the 2-category of (small) categories. However, we have used the terminology 2-universal more broadly than \cite{Fio}, as they work exclusively with strict 2-categories. Given that every 2-category is naturally equivalent to a strict one, this is not a cause for concern.
\end{Remark}

\begin{Lemma}\label{lem:1aroubienvelopeexist}
Given any (1-)category $\mathcal{C}$, its Karoubi envelope exists.
\end{Lemma}
\begin{proof}
The existence of a functor $\iota_{\mathcal{C}}:\mathcal{C}\rightarrow Kar(\mathcal{C})$ is standard. Verifying the 2-universal property follows from the universal property of splittings of idempotents, and the fact that splittings of idempotents are preserved by all functors.
\end{proof}

The following result is part of the folklore.

\begin{Proposition}\label{prop:Karfunctor}
The Karoubi envelope defines a 2-functor $$Kar:\mathbf{2Cat}\rightarrow \mathbf{2Cat},$$ and a 2-natural transformation $$\iota: Id \Rightarrow Kar.$$ Further, if $\mathcal{C}$ is an idempotent complete category, $\iota_{\mathcal{C}}$ is an equivalence of categories.
\end{Proposition}

\begin{Remark}\label{rem:3universalKarfunctor}
To be precise, the definitions of $Kar$ and $\iota$ involve choices. However, using remark \ref{rem:2univbiuniv}, and theorem 9.16 of \cite{Fio}, we obtain the following result, which shows that these choices are all essentially unique.
\end{Remark}

\begin{Proposition}\label{prop:Karleft2adjoint}
The 2-functor $Kar:\mathbf{2Cat}\rightarrow \mathbf{2Cat}^{ic}$ is a left 2-adjoint for the forgetful 2-functor $U$, with unit given by $\iota$.
\end{Proposition}

It is useful to know when the Karoubi envelope of a functor is an equivalence. The following characterization is standard.

\begin{Lemma}\label{lem:Karequivalence}
Let $F:\mathcal{C}\rightarrow \mathcal{D}$ be a functor. Then, $Kar(F)$ is an equivalence if and only if $F$ is fully faithful, and every object of $\mathcal{D}$ can be expressed as the splitting of an idempotent supported on an object in the image of $F$.
\end{Lemma}

Let us presently indicate the changes that are necessary to accommodate $R$-linear categories, where $R$ is a fixed commutative ring. In this context, the appropriate notion of completion is called the Cauchy completion, and involves both the splitting of idempotents and the existence of finite direct sums. The Cauchy completion is given by first taking the direct sum completion, and then the Karoubi envelope. Given an $R$-linear (1-)category $\mathcal{C}$, we denote its Cauchy completion by $\kappa_{\mathcal{C}}:\mathcal{C}\rightarrow Cau(\mathcal{C})$. The $R$-linear functor $\kappa_{\mathcal{C}}$ has a 2-universal property similar to the one of the Karoubi envelope.

\begin{Lemma}\label{lem:Cauchy1completion}
Every $R$-linear category has a Cauchy completion.
\end{Lemma}
\begin{proof}
The existence of the direct sum completion is standard. Combining this observation with lemma \ref{lem:1aroubienvelopeexist}, and the fact that the Karoubi envelope of an $R$-linear category that has finite direct sums is $R$-linear and has finite direct sums yields the result.
\end{proof}

Let us denote by $\mathbf{2Cat}_{R}$ the $R$-linear 2-category of $R$-linear categories, and by $\mathbf{2Cat}_{R}^{cc}$ the full sub-2-category on the Cauchy complete $R$-linear categories. Using the the 2-universal property of the Cauchy completion, and \cite[Thm. 9.16]{Fio}, we have:

\begin{Proposition} \label{lem:directsumcompletion2adjoint}
Cauchy completion defines a $R$-linear 2-functor $$Cau:\mathbf{2Cat}_{R}\rightarrow \mathbf{2Cat}_{R}^{cc}$$ that is left 2-adjoint to the forgetful 2-functor $V:\mathbf{2Cat}_{R}^{cc}\rightarrow \mathbf{2Cat}_{R}$, with unit $\kappa$.
\end{Proposition}

It is not hard to deduce a characterization of those $R$-linear functors that become equivalences upon application of $Cau$ from lemma \ref{lem:Karequivalence}.

\begin{Lemma}\label{lem:Cauchyequivalence}
Let $F:\mathcal{C}\rightarrow \mathcal{D}$ be a $R$-linear functor of $R$-linear categories. Then, $Cau(F)$ is an equivalence if and only if $F$ is fully faithful and every object of $\mathcal{D}$ is the splitting of an idempotent supported on a finite direct sum of objects in the image of $F$.
\end{Lemma}

\subsection{Completions of 2-Categories}\label{sec:completion2categories}
Let us now turn our attention to the 2-categorical story. From now on, we shall assume that all 2-categories under consideration are locally idempotent complete, i.e. their $Hom$-categories are idempotent complete. The goal of this section is to prove proposition \ref{prop:exsitenceKaroubienvelope}. Before doing so, we need to unfold the following definition from \cite{GJF}.

\begin{Definition}\label{def:2condensationbimodule}
Let $\mathfrak{C}$ be a 2-category, $(A_1,e_1,\mu_1,\delta_1)$ and $(A_2,e_2,\mu_2,\delta_2)$ be 2-condensation monads. An $(A_2,A_1)$-bimodule consists of:

\begin{enumerate}
    \item A 1-morphism $b:A_1\rightarrow A_2$;
    \item Two 2-morphisms $\nu^r: b\circ e_1\Rightarrow b$ and $\beta^r: b\Rightarrow b\circ e_1$;
    \item Two 2-morphisms $\nu^l: e_2\circ b\Rightarrow b$ and $\beta^l: b\Rightarrow e_2\circ b$
\end{enumerate}

such that

\begin{enumerate}
    \item The pair $(f,\nu^r, \nu^l)$ is an associative $(e_2,e_1)$-bimodule;
    \item The pair $(f,\beta^r, \beta^l)$ is a coassociative $(e_2,e_1)$-bicomodule;
    \item The 2-morphisms $\nu^r$ and $\beta^l$ commute, and the 2-morphisms $\nu^l$ and $\beta^r$ commute;
    \item The Frobenius relations for modules hold for $(b,\nu^r,\beta^r)$ and $(b,\nu^l,\beta^l)$;
    \item We have $\beta^r \cdot \nu^r = Id_b$ and $\beta^l \cdot \nu^l = Id_b$.
\end{enumerate}
\end{Definition}

\begin{Notation}
When working with bimodules as defined above in a strict 2-category, it is extremely convenient to use the graphical notation given in example 3.0.1 of \cite{GJF}. Further, as we shall exclusively read our diagrams from right to left, and from bottom to top, we suppress the arrows from the notation. Using these conventions, the Frobenius relations for $(b,\nu^r,\beta^r)$ can be expressed as the following equalities:

$$\begin{tikzpicture}[xscale=0.5, yscale=0.4]
\draw[line width=0.5mm] (0,0) node[anchor=north] {$b_1$}  -- (0,3) node[anchor=south] {$b_1$};
\draw (2,0) node[anchor=north] {$e_1$} -- (0,2);
\draw (1,1) -- (2,3)node[anchor=south] {$e_1$};

\draw (3,1.5) node {$=$};

\draw[line width=0.5mm] (4,0) node[anchor=north] {$b_1$}  -- (4,3) node[anchor=south] {$b_1$};
\draw (5,0) node[anchor=north west] {$e_1$} -- (4,1);
\draw (4,2) -- (5,3)node[anchor=south west] {$e_1$};

\draw (6,1.5) node {$=$};

\draw[line width=0.5mm] (7,0) node[anchor=north] {$b_1$}  -- (7,3) node[anchor=south] {$b_1$};
\draw (9,0) node[anchor=north] {$e_1$} -- (8,2);
\draw (7,1) -- (9,3)node[anchor=south] {$e_1$};

\draw (9.5,1.5) node {.};
\end{tikzpicture}$$

\end{Notation}

\renewcommand*{\proofname}{Proof of proposition \ref{prop:exsitenceKaroubienvelope}}

\begin{proof}
Without loss of generality, we may assume that $\mathfrak{C}$ is a strict 2-category. The condensation complete 2-category $Kar(\mathfrak{C})$ and the 2-functor $\iota_{\mathcal{C}}$ are defined in theorem 2.3.10 of \cite{GJF}. In particular, recall that $Kar(\mathfrak{C})$ is the 2-category of 2-condensation monads in $\mathfrak{C}$, with 1-morphisms given by bimodules, and 2-morphisms by bimodule maps. The 2-functor $\iota_{\mathfrak{C}}$ sends an object $C$ of $\mathfrak{C}$ to the trivial 2-condensation monad $(C, Id_C, Id_{Id_C}, Id_{Id_C})$ on $C$.

Let $F:\mathfrak{C}\rightarrow \mathfrak{A}$ be a 2-functor to a condensation complete 2-category $\mathfrak{A}$. Without loss of generality, we may assume that $\mathfrak{A}$ is strict. Given 1-morphisms $c:C\rightarrow D$, and $d:D\rightarrow E$ in $\mathfrak{C}$, we write the coherence 2-isomorphism witnessing that $F$ respects the composition by $F_{d,c}:F(d) \circ F(c)\cong F(d\circ c)$. Further, given an object $C$ of $\mathfrak{C}$, we denote the unitor by $F_C: Id_{F(C)}\cong F(Id_C)$.

We now define a 2-functor $F':Kar(\mathfrak{C})\rightarrow \mathfrak{A}$. On objects, it is defined using theorem 2.3.1 of \cite{GJF}, which says that the 2-category of extensions to a 2-condensation of a given 2-condensation monad is either empty or contractible. Because $\mathfrak{A}$ has all condensates, they are all contractible. Explicitly, given 2-condensation monads $(A_i,e_i,\mu_i,\delta_i)$ in $\mathfrak{C}$ for $i=1,2,3,4$, their images under $F$ are the 2-condensation monads $$F(A_i,e_i,\mu_i,\delta_i):=(F(A_i),F(e_i),F(\mu_i)\cdot F_{e_i,e_i},F_{e_i,e_i}^{-1}\cdot F(\delta_i))$$ in $\mathfrak{A}$. For each $i$, we fix a 2-condensation $(F(A_i),B_i,f_i,g_i,\phi_i,\gamma_i)$ in $\mathfrak{A}$ together with a 2-isomorphism $\theta_i:g_i\circ f_i\cong e_i$ providing an extension of $F(A_i,e_i,\mu_i,\nu_i)$ in $\mathfrak{A}$ to 2-condensations, and we set $F'(A_i) := B_i$.

In order to define $F'$ on 1-morphisms, we use the fact that $\mathfrak{A}$ is locally idempotent complete. For $j=1,2,3$, let $(b_j, \nu^r_j, \beta^r_j, \nu^l_j, \beta^l_j)$ be an $(A_{j+1}, A_j)$ bimodules. The 1-morphism $f_2\circ F(b_1)\circ g_1$ supports the idempotent 2-morphism specified by:

$$
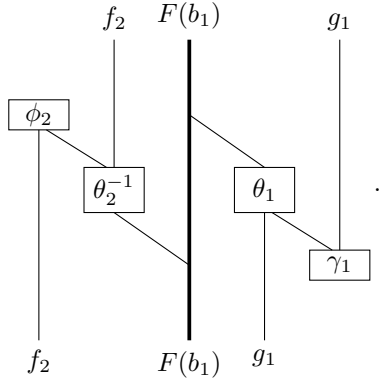
\begin{tikzpicture}[xscale=1, yscale=1]
\draw (0,0) node[anchor=north] {$f_2$} -- (0,2.8);

\draw (0,3) node{$\phi_2$};
\draw (-0.4,2.8) rectangle (0.4,3.2);

\draw (1,2) node{$\theta^{-1}_2$};
\draw (0.6,1.7) rectangle (1.4,2.3);

\draw (0.9,2.3) -- (0.1,2.8);
\draw (2,1) -- (1, 1.7);
\draw (1,2.3) -- (1,4)node[anchor=south] {$f_2$};

\draw[line width=0.5mm] (2,0) node[anchor=north] {$F(b_1)$}  -- (2,4) node[anchor=south] {$F(b_1)$};

\draw (3.9,1.2) -- (3.1,1.7);
\draw (3,0) node[anchor=north] {$g_1$} -- (3, 1.7);
\draw (3,2.3) -- (2,3);

\draw (3,2) node{$\theta_1$};
\draw (2.6,1.7) rectangle (3.4,2.3);

\draw (4,1) node{$\gamma_1$};
\draw (3.6,0.8) rectangle (4.4,1.2);

\draw (4,1.2)  -- (4,4)node[anchor=south] {$g_1$};

\draw (4.5,2) node{.};

\end{tikzpicture}$$

\noindent The 1-morphism $F'(b_1)$ is given by a choice of splitting for this idempotent in $Hom_{\mathfrak{A}}(B_1,B_2)$. The same procedure is used to define $F'(b_2)$ and $F'(b_3)$. The value of $F'$ on 2-morphisms is defined similarly using the universal property of splittings of idempotents. In order to simplify notations, from now on, we shall omit the subscript in the expressions labelling the boxes.

It remains to prove that $F'$ is a 2-functor and that there exists a 2-natural equivalence $F'\circ\iota_{\mathfrak{C}}\Rightarrow F$. Recall that the composite of $b_1$ and $b_2$ in $Kar(\mathfrak{C})$, denoted by $b_2\otimes_{A_2}b_1$, is obtained by choosing a splitting for the idempotent

$$
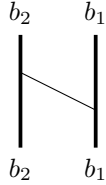
\begin{tikzpicture}[xscale=0.5, yscale=0.5]

\draw[line width=0.5mm] (0,0) node[anchor=north] {$b_{2}$}  -- (0,3) node[anchor=south] {$b_2$};

\draw (2,1) -- (0,2);

\draw[line width=0.5mm] (2,0) node[anchor=north] {$b_1$}  -- (2,3) node[anchor=south] {$b_1$};

\end{tikzpicture}$$

\noindent supported on $b_2\circ b_1$ (see \cite{GJF} for details). Thus, we see that both $F'(b_2)\circ F'(b_1)$ and $F'(b_2\otimes_{A_2}b_1)$ are defined as splittings of idempotents. Using the fact that such splittings of idempotent are universal, we find that $F'(b_2)\circ F'(b_1)$ is a splitting of the idempotent supported on $f_3\circ F(b_2)\circ g_2\circ f_2\circ F(b_1)\circ g_1$ given by:

$$
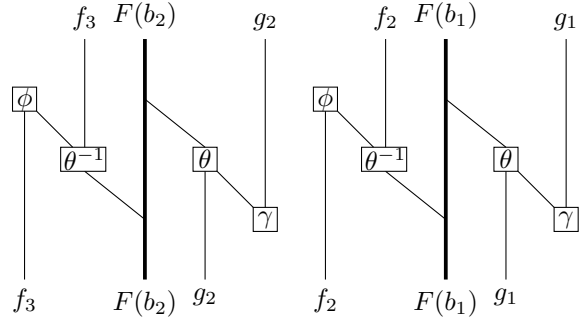
\begin{tikzpicture}[xscale=0.8, yscale=0.8]
\draw (0,0) node[anchor=north] {$f_3$} -- (0,2.8);

\draw (-0.2,2.8) rectangle (0.2,3.2);
\draw (0,3) node{$\phi$};

\draw (0.6,1.8) rectangle (1.35,2.2);
\draw (1,2) node{$\theta^{-1}$};

\draw (0.8,2.2) -- (0.2,2.8);
\draw (2,1) -- (1, 1.8);
\draw (1,2.2) -- (1,4)node[anchor=south] {$f_3$};

\draw[line width=0.5mm] (2,0) node[anchor=north] {$F(b_2)$}  -- (2,4) node[anchor=south] {$F(b_2)$};

\draw (3.8,1.2) -- (3.2,1.8);
\draw (3,0) node[anchor=north] {$g_2$} -- (3, 1.8);
\draw (3,2.2) -- (2,3);

\draw (2.8,1.8) rectangle (3.2,2.2);
\draw (3,2) node{$\theta$};

\draw (3.8,0.8) rectangle (4.2,1.2);
\draw (4,1) node{$\gamma$};

\draw (4,1.2)  -- (4,4)node[anchor=south] {$g_2$};

\draw (5,0) node[anchor=north] {$f_2$} -- (5,2.8);

\draw (4.8,2.8) rectangle (5.2,3.2);
\draw (5,3) node{$\phi$};

\draw (5.6,1.8) rectangle (6.35,2.2);
\draw (6,2) node{$\theta^{-1}$};

\draw (5.8,2.2) -- (5.2,2.8);
\draw (7,1) -- (6, 1.8);
\draw (6,2.2) -- (6,4)node[anchor=south] {$f_2$};

\draw[line width=0.5mm] (7,0) node[anchor=north] {$F(b_1)$}  -- (7,4) node[anchor=south] {$F(b_1)$};

\draw (8.8,1.2) -- (8.2,1.8);
\draw (8,0) node[anchor=north] {$g_1$} -- (8, 1.8);
\draw (8,2.2) -- (7,3);

\draw (7.8,1.8) rectangle (8.2,2.2);
\draw (8,2) node{$\theta$};

\draw (8.8,0.8) rectangle (9.2,1.2);
\draw (9,1) node{$\gamma$};

\draw (9,1.2)  -- (9,4)node[anchor=south] {$g_1$};

\draw (9.5,2) node{.};

\end{tikzpicture}$$

\noindent Likewise, one find that $F'(b_2\otimes_{A_2}b_1)$ is a splitting of the idempotent supported on $f_3\circ F(b_2)\circ F(b_1)\circ g_1$ specified by:

$$
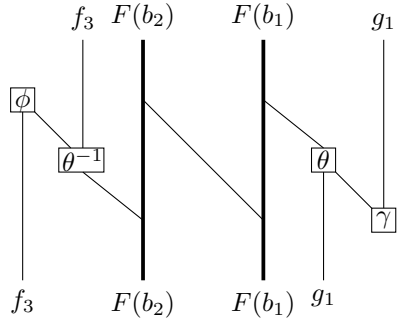
\begin{tikzpicture}[xscale=0.8, yscale=0.8]
\draw (0,0) node[anchor=north] {$f_3$} -- (0,2.8);

\draw (-0.2,2.8) rectangle (0.2,3.2);
\draw (0,3) node{$\phi$};

\draw (0.6,1.8) rectangle (1.35,2.2);
\draw (1,2) node{$\theta^{-1}$};

\draw (0.8,2.2) -- (0.2,2.8);
\draw (2,1) -- (1, 1.8);
\draw (1,2.2) -- (1,4)node[anchor=south] {$f_3$};

\draw[line width=0.5mm] (2,0) node[anchor=north] {$F(b_2)$}  -- (2,4) node[anchor=south] {$F(b_2)$};

\draw (4,1) -- (2,3);

\draw[line width=0.5mm] (4,0) node[anchor=north] {$F(b_1)$}  -- (4,4) node[anchor=south] {$F(b_1)$};

\draw (5.8,1.2) -- (5.2,1.8);
\draw (5,0) node[anchor=north] {$g_1$} -- (5, 1.8);
\draw (5,2.2) -- (4,3);

\draw (4.8,1.8) rectangle (5.2,2.2);
\draw (5,2) node{$\theta$};

\draw (5.8,0.8) rectangle (6.2,1.2);
\draw (6,1) node{$\gamma$};

\draw (6,1.2)  -- (6,4)node[anchor=south] {$g_1$};

\draw (6.5,2) node{.};

\end{tikzpicture}$$

\noindent Now, observe that the 2-morphisms with source $f_3\circ F(b_2)\circ F(b_1)\circ g_1$ and target $f_3\circ F(b_2)\circ g_2\circ f_2\circ F(b_1)\circ g_1$ given by 

$$
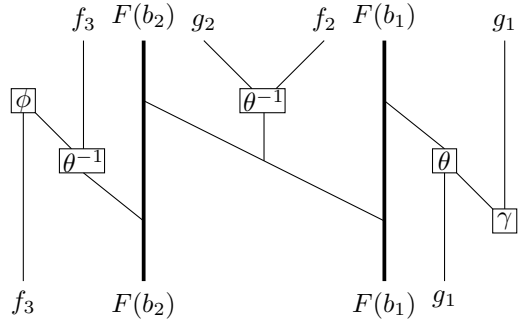
\begin{tikzpicture}[xscale=0.8, yscale=0.8]
\draw (0,0) node[anchor=north] {$f_3$} -- (0,2.8);

\draw (-0.2,2.8) rectangle (0.2,3.2);
\draw (0,3) node{$\phi$};

\draw (0.6,1.8) rectangle (1.35,2.2);
\draw (1,2) node{$\theta^{-1}$};

\draw (0.8,2.2) -- (0.2,2.8);
\draw (2,1) -- (1, 1.8);
\draw (1,2.2) -- (1,4)node[anchor=south] {$f_3$};

\draw[line width=0.5mm] (2,0) node[anchor=north] {$F(b_2)$}  -- (2,4) node[anchor=south] {$F(b_2)$};

\draw (6,1) -- (2,3);
\draw (4,2) -- (4,2.8);

\draw (3.6,2.8) rectangle (4.35,3.2);
\draw (4,3) node{$\theta^{-1}$};

\draw (3.8,3.2) -- (3,4)node[anchor=south] {$g_2$};
\draw (4.2,3.2) -- (5,4)node[anchor=south] {$f_2$};

\draw[line width=0.5mm] (6,0) node[anchor=north] {$F(b_1)$}  -- (6,4) node[anchor=south] {$F(b_1)$};

\draw (7.8,1.2) -- (7.2,1.8);
\draw (7,0) node[anchor=north] {$g_1$} -- (7, 1.8);
\draw (7,2.2) -- (6,3);

\draw (6.8,1.8) rectangle (7.2,2.2);
\draw (7,2) node{$\theta$};

\draw (7.8,0.8) rectangle (8.2,1.2);
\draw (8,1) node{$\gamma$};

\draw (8,1.2)  -- (8,4)node[anchor=south] {$g_1$};

\end{tikzpicture}$$

\noindent induces an equivalence of idempotents between the two idempotents defined above. Namely, its inverse is given by the 2-morphism 

$$
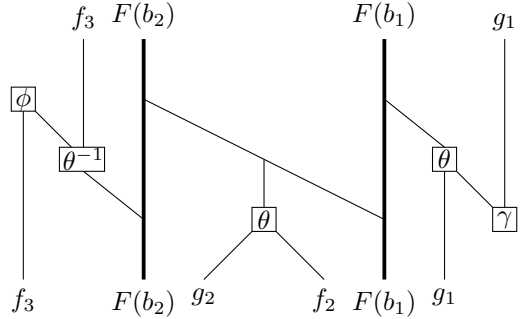
\begin{tikzpicture}[xscale=0.8, yscale=0.8]
\draw (0,0) node[anchor=north] {$f_3$} -- (0,2.8);

\draw (-0.2,2.8) rectangle (0.2,3.2);
\draw (0,3) node{$\phi$};

\draw (0.6,1.8) rectangle (1.35,2.2);
\draw (1,2) node{$\theta^{-1}$};

\draw (0.8,2.2) -- (0.2,2.8);
\draw (2,1) -- (1, 1.8);
\draw (1,2.2) -- (1,4)node[anchor=south] {$f_3$};

\draw[line width=0.5mm] (2,0) node[anchor=north] {$F(b_2)$}  -- (2,4) node[anchor=south] {$F(b_2)$};

\draw (6,1) -- (2,3);
\draw (4,1.2) -- (4,2);

\draw (3.8,0.8) rectangle (4.2,1.2);
\draw (4,1) node{$\theta$};

\draw (3,0) node[anchor=north] {$g_2$} -- (3.8, 0.8);
\draw (5,0) node[anchor=north] {$f_2$} -- (4.2, 0.8);

\draw[line width=0.5mm] (6,0) node[anchor=north] {$F(b_1)$}  -- (6,4) node[anchor=south] {$F(b_1)$};

\draw (7.8,1.2) -- (7.2,1.8);
\draw (7,0) node[anchor=north] {$g_1$} -- (7, 1.8);
\draw (7,2.2) -- (6,3);

\draw (6.8,1.8) rectangle (7.2,2.2);
\draw (7,2) node{$\theta$};

\draw (7.8,0.8) rectangle (8.2,1.2);
\draw (8,1) node{$\gamma$};

\draw (8,1.2)  -- (8,4)node[anchor=south] {$g_1$};

\end{tikzpicture}$$

\noindent from $f_3\circ F(b_2)\circ g_2\circ f_2\circ F(b_1)\circ g_1$ to $f_3\circ F(b_2)\circ F(b_1)\circ g_1$. Thus, there exists a unique 2-isomorphism $F'_{b_2,b_1}:F'(b_2)\circ F'(b_1)\cong F'(b_2\otimes_{A_2}b_1)$. Further, because of the uniqueness of the isomorphism induced by an equivalence of idempotents, we deduce that $$F'_{b_3\otimes_{A_3} b_2,b_1}\cdot (F'_{b_3,b_2}\circ F'(b_1)) = F'_{b_3,b_2\otimes_{A_2}b_1}\cdot (F'(b_3)\circ F'_{b_2,b_1}).$$ This proves that $F'$ preserves the composition of 1-morphisms. A similar argument shows that $F'$ preserves units, whence we find that $F'$ is indeed a 2-functor.

We still have to construct a 2-natural equivalence $u:F'\circ\iota_{\mathfrak{C}}\Rightarrow F$. Given $i=1,2$, observe that on the objects $C_i$ of $\mathfrak{C}$, $F' \iota_{\mathfrak{C}}(C_i)$ is given by an extension of the image under $F$ of the 2-condensation monad $(C_i, Id_{C_i}, Id_{Id_{C_i}}, Id_{Id_{C_i}})$. Note that the 2-condensation $(F({C_i}), F({C_i}), Id_{F(C_i)}, Id_{F(C_i)}, Id_{Id_{F(C_i)}},Id_{Id_{F(C_i)}})$ together with the 2-isomorphism $F_{C_i}$ is such an extension. But, we have already picked an extension above, namely $(F(C_i), B_i,f_i,g_i,\phi_i,\gamma_i)$. By the proof of theorem 2.3.2 of \cite{GJF}, we get that the 1-morphism $g_i:B_i\rightarrow F(C_i)$ is an equivalence $F' \iota_{\mathfrak{C}}({C_i})\rightarrow F({C_i})$ in $\mathfrak{A}$ for $i=1,2$. Hence, it is sufficient to prove that these $g_i$ can be assembled into a 2-natural transformation $u$.

Let $c:C_1\rightarrow C_2$, be a 1-morphism in $\mathfrak{C}$. Abusing notation slightly by writing $C_i$, for the 2-condensation monad $(C_i, Id_{C_i}, Id_{Id_{C_i}}, Id_{Id_{C_i}})$ with $i=1,2$, it is clear that we can endow $c$ with a $(C_2,C_1)$-bimodule structure using identity morphisms. This is precisely the value of $\iota_{\mathfrak{C}}(c)$. Then, by definition, one gets that $g_2\circ(F'\iota_{\mathfrak{C}})(c)$ is a splitting of the idempotent supported on $g_2\circ f_2\circ F(c)\circ g_1$ given by:

$$
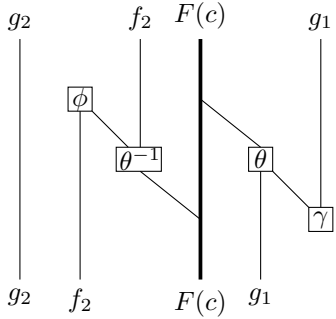
\begin{tikzpicture}[xscale=0.8, yscale=0.8]
\draw (-1,0) node[anchor=north] {$g_2$}  -- (-1,4) node[anchor=south] {$g_2$};

\draw (0,0) node[anchor=north] {$f_2$} -- (0,2.8);

\draw (-0.2,2.8) rectangle (0.2,3.2);
\draw (0,3) node{$\phi$};

\draw (0.6,1.8) rectangle (1.35,2.2);
\draw (1,2) node{$\theta^{-1}$};

\draw (0.8,2.2) -- (0.2,2.8);
\draw (2,1) -- (1, 1.8);
\draw (1,2.2) -- (1,4)node[anchor=south] {$f_2$};

\draw[line width=0.5mm] (2,0) node[anchor=north] {$F(c)$}  -- (2,4) node[anchor=south] {$F(c)$};

\draw (3.8,1.2) -- (3.2,1.8);
\draw (3,0) node[anchor=north] {$g_1$} -- (3, 1.8);
\draw (3,2.2) -- (2,3);

\draw (2.8,1.8) rectangle (3.2,2.2);
\draw (3,2) node{$\theta$};

\draw (3.8,0.8) rectangle (4.2,1.2);
\draw (4,1) node{$\gamma$};

\draw (4,1.2)  -- (4,4)node[anchor=south] {$g_1$};

\draw (4.5,2) node{.};

\end{tikzpicture}$$

\noindent We claim that $F(c_1)\circ g_1$ is also a splitting of the above idempotent. This can be seen by using the two 2-morphisms given by the two diagrams below:

$$
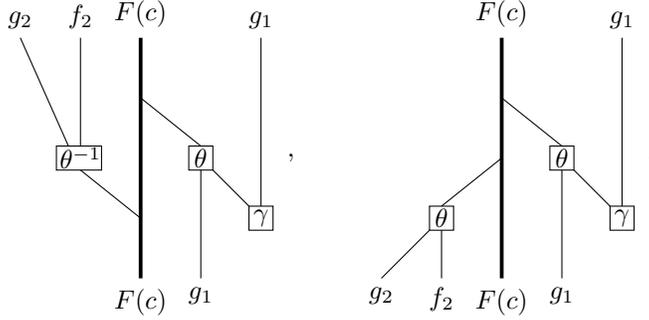
\begin{tikzpicture}[xscale=0.8, yscale=0.8]

\draw (0.6,1.8) rectangle (1.35,2.2);
\draw (1,2) node{$\theta^{-1}$};

\draw (0.8,2.2) -- (0,4)node[anchor=south] {$g_2$};
\draw (2,1) -- (1, 1.8);
\draw (1,2.2) -- (1,4)node[anchor=south] {$f_2$};

\draw[line width=0.5mm] (2,0) node[anchor=north] {$F(c)$}  -- (2,4) node[anchor=south] {$F(c)$};

\draw (3.8,1.2) -- (3.2,1.8);
\draw (3,0) node[anchor=north] {$g_1$} -- (3, 1.8);
\draw (3,2.2) -- (2,3);

\draw (2.8,1.8) rectangle (3.2,2.2);
\draw (3,2) node{$\theta$};

\draw (3.8,0.8) rectangle (4.2,1.2);
\draw (4,1) node{$\gamma$};

\draw (4,1.2)  -- (4,4)node[anchor=south] {$g_1$};

\draw (4.5,2) node{,};

\draw (6.8,0.8) rectangle (7.2,1.2);
\draw (7,1) node{$\theta$};

\draw (6,0)node[anchor=north] {$g_2$} -- (6.8,0.8);
\draw (7,1.2) -- (8,2);
\draw (7,0)node[anchor=north] {$f_2$} -- (7,0.8);

\draw[line width=0.5mm] (8,0) node[anchor=north] {$F(c)$}  -- (8,4) node[anchor=south] {$F(c)$};

\draw (9.8,1.2) -- (9.2,1.8);
\draw (9,0) node[anchor=north] {$g_1$} -- (9, 1.8);
\draw (9,2.2) -- (8,3);

\draw (8.8,1.8) rectangle (9.2,2.2);
\draw (9,2) node{$\theta$};

\draw (9.8,0.8) rectangle (10.2,1.2);
\draw (10,1) node{$\gamma$};

\draw (10,1.2)  -- (10,4)node[anchor=south] {$g_1$};

\draw (10.5,2) node{.};

\end{tikzpicture}$$

\noindent Hence, there is a unique 2-isomorphism $$u_{c}:g_2\circ(F'\iota_{\mathfrak{C}})(c)\cong F(c) \circ g_1.$$ Using the uniqueness of the 2-isomorphism comparing two splittings of a given idempotent, it follows that the 1-morphisms $g_1$, $g_2$ together with the 2-isomorphisms $u_{c}$ for varying $c$ define a 2-natural transformation $u$. As the 1-morphisms $g_1$, and $g_2$ are equivalences, we get a 2-natural equivalence $u:F'\circ \iota_{\mathfrak{C}}\Rightarrow F$.

The second and third point of definition \ref{def:KarEnv} follow using analogous arguments.
\end{proof}

\renewcommand*{\proofname}{Proof}

\begin{Remark}
As a cageorification of proposition \ref{prop:Karfunctor}, it is possible to prove that $Kar$ defines a 3-functor from the 3-category of 2-categories to itself, and that $\iota$ defines a 3-natural transformation $\iota:Id\Rightarrow Kar$. For an explicit proof, albeit in a somewhat different context, we refer the reader to \cite{CP}. Further, we expect that the results of section 9 of \cite{Fio} can be categorified. Together with proposition \ref{prop:exsitenceKaroubienvelope}, this would yield a proof of the fact that $Kar$, viewed as a 3-functor to the 3-category of condensation complete 2-categories is a left 3-adjoint for the forgetful 3-functor.
\end{Remark}

It will be important to have a categorified version of lemma \ref{lem:Karequivalence}. Namely, we want to characterize those 2-functors that become 2-equivalences upon taking the Karoubi envelope.

\begin{Lemma}\label{lem:2Karoubiequivalence}
Let $F:\mathfrak{C}\rightarrow \mathfrak{D}$ be a 2-functor. Then, $Kar(F)$ is an equivalence of 2-categories if and only if $F$ is fully faithful (it induces a equivalences on $Hom$-categories), and every object of $\mathfrak{D}$ is the splitting of a condensation monad on an object in the image of $F$.
\end{Lemma}

\begin{proof}
The backward implication follows from the fact that the 2-functors $\iota_{\mathfrak{C}}$ and $\iota_{\mathfrak{D}}$ are fully faithful (see theorem 2.3.10 of \cite{GJF}). Conversely, let us assume that $F$ has the above properties. Then, it follows from the construction, and the fact that $F$ is fully faithful, that $Kar(F)$ is fully faithful. It remains to prove the essential surjectivity. In spirit, the argument is similar to the proof that the Karoubi envelope of a 2-category is condensation complete (see proposition A.5.3 of \cite{DR}). Let $(D,e, \mu, \delta)$ be a 2-condensation monad in $\mathfrak{D}$. By hypothesis, there exists a 2-condensation $$(F(C),D, f,g, \phi,\gamma)$$ in $\mathfrak{D}$ for some $C$ in $\mathfrak{C}$. We can compose the 2-condensation monad on $D$ with this 2-condensation, and produce a 2-condensation monad $$(F(C), g\circ e\circ f, (f\circ \mu\circ g)\cdot (g\circ e\circ \phi\circ e\circ f), (g\circ e\circ \gamma\circ e\circ f)\cdot(g\circ \delta\circ f)).$$ Now, note that these 2-condensation monads are equivalent as objects of $Kar(\mathfrak{D})$ via the canonical 2-condensation bimodules given by $f$, and $g$. As $F$ is fully faithful, the 2-condensation monad $(F(C), g\circ e\circ f)$ is the image of a 2-condensation monad in $\mathfrak{C}$. This proves the result.
\end{proof}

Let $R$ be a fixed commutative ring. Recall that an $R$-linear locally additive (and locally idempotent complete) 2-category is Cauchy complete provided it has all condensates and finite direct sums. Theorem 4.3.5 of \cite{GJF} shows that for $R$-linear locally additive 2-categories and $R$-linear 2-functors, the natural notion of completion is the Cauchy completion.

The Cauchy completion of a locally additive $R$-linear 2-category $\mathfrak{C}$ can be constructed as follows: First, take the direct sum completion of $\mathfrak{C}$, denoted by $Mat(\mathfrak{C})$. This is the locally additive 2-category whose objects are finite vectors with values in $\mathfrak{C}$, and whose morphisms are finite matrices of morphisms in $\mathfrak{C}$. Observe that the canonical 2-functor $\mathfrak{C}\rightarrow Mat(\mathfrak{C})$ has an obvious 3-universal property. Second, take the Karoubi envelope. We get an $R$-linear 2-functor $\mathfrak{C}\rightarrow Kar(Mat(\mathfrak{\mathfrak{C}}))$. Noting that the Karoubi envelope of a locally additive $R$-linear 2-category with finite direct sums has finite direct sums shows that the target 2-category is indeed Cauchy complete.

Given a locally additive $R$-linear 2-category $\mathfrak{C}$, we denote the Cauchy completion constructed above by $\kappa_\mathfrak{C}:\mathfrak{C}\rightarrow Cau(\mathfrak{C})$. It enjoys a 3-universal property analogous to the one of the Karoubi envelope. Namely, it is enough to modify definition \ref{def:KarEnv} by replacing $\iota_{\mathfrak{C}}$ by $\kappa_\mathfrak{C}:\mathfrak{C}\rightarrow Cau(\mathfrak{C})$, requiring that all 2-categories and 2-functors be $R$-linear, and that $\mathfrak{A}$ be Cauchy complete.

Finally, in a subsequent article, we will need to know when the Cauchy completion of an $R$-linear 2-functor is an equivalence. The answer is given by the next lemma, which is a consequence of lemma \ref{lem:2Karoubiequivalence}.

\begin{Lemma}\label{lem:2Cauchyequivalence}
Let $F:\mathfrak{C}\rightarrow \mathfrak{D}$ be an $R$-linear 2-functor of locally additive $R$-linear 2-categories. Then, $Cau(F)$ is an equivalence of 2-categories if and only if $F$ is fully faithful, and every object of $\mathfrak{D}$ is the splitting of a 2-condensation monad supported on a direct sum of objects in the image of $F$.
\end{Lemma}

\subsubsection*{Acknowledgments}

I would like to thank Christopher Douglas for helpful conversations related to the content of this article.

\bibliography{bibliography.bib}

\end{document}